\documentclass[12pt]{article}

\usepackage{amssymb, amscd}
\usepackage{amsmath}
\usepackage{amsthm}
\usepackage{epsfig}
\usepackage{psfig}
\usepackage{xcolor}

\newtheorem{dfn}{Definition}[section]

\newtheorem{rem}[dfn]{Remark}

\newtheorem{definition}[dfn]{Definition}

\newtheorem{thm}[dfn]{Theorem}
\newtheorem{lem}[dfn]{Lemma}
\newtheorem{assumption}[dfn]{Assumption}
\newtheorem{notation}[dfn]{Notation}

\newtheorem{prop}[dfn]{Proposition}
\newtheorem{cor}[dfn]{Corollary}

\newtheorem{example}[dfn]{Example}

\def\proof{\par\medskip\noindent{\it Proof: }}

\newcommand{\Isom}{\operatorname{Isom}}

\newcommand{\diam}{\operatorname{diam}}
\newcommand{\length}{\operatorname{length}}
\newcommand{\dist}{\operatorname{dist}}

\def\D{\partial}
\def\la{\lambda}
\def\R{{\mathbb R}}
\def\C{{\mathbb C}}
\def\N{{\mathbb N}}

\def\e{{\mathcal E}}
\def\f{{\mathcal F}}

\def\n{{\mathcal N}}
\def\v{{\mathcal V}}

\def\Ga{\Gamma}
\def\eps{\epsilon}
\def\al{\alpha}
\def\be{\beta}
\def\ga{\gamma}

\def\ol{\overline}

\def\acts{\curvearrowright}
\def\<{\langle}
\def\>{\rangle}

\begin{document}

\title{Energy of harmonic functions and Gromov's proof of Stallings' theorem}
\author{M.~Kapovich}
\date{\today}

\maketitle

\begin{abstract}
We provide the details for Gromov's proof of Stallings' theorem on groups
with infinitely many ends using harmonic functions.
The main technical result of the paper is a compactness theorem
for a certain family of harmonic functions.
\end{abstract}

\section{Introduction}

In his essay \cite[Pages 228--230]{Gromov}, Gromov gave a proof of the Stallings' theorem
on  groups with infinitely many ends using harmonic functions.
The goal of this paper is to
provide the details for Gromov's arguments. The main bulk of the paper
is devoted to the proof of a compactness theorem for a certain family
of harmonic functions. The corresponding statements are contained in
Steps 2 and 4 of Gromov's argument. The rest of our proof closely follows
Gromov's. 

Let $M$ be a complete Riemannian manifold of bounded geometry, which has
infinitely many ends. Suppose that there exists a number $R$ such that
every point in $M$ belongs to an $R$--{\em neck}, i.e., an $R$-ball which
separates $M$ into at least three unbounded components.
(This property is immediate if $M$ admits a cocompact isometric group action.)

Let $\ol{M}:= M\cup Ends(M)$ denote the compactification of $M$ by its space of ends.
Given a continuous function $\chi: Ends(M)\to \{0, 1\}$,
let
$$
h=h_\chi: \ol{M} \to [0, 1]
$$
denote the continuous extension of $\chi$, so that $h|M$ is
harmonic. Let $H(M)$ denote the space of harmonic functions
$$
\{h=h_{\chi}, \chi: Ends(M)\to \{0, 1\}
\hbox{~~is nonconstant} \}.
$$
We give $H(M)$ the topology of uniform convergence on compacts in $M$. Let
$E: H(M)\to \R_+=[0, \infty)$ denote the energy functional.

\begin{definition}
Given the manifold $M$, define its {\em energy gap} $e(M)$ as
$$
e(M):=\inf \{ E(h): h\in H(M) \}.
$$
\end{definition}

If $M$ admits an isometric group action $G\acts M$, then $G$ acts on $H(M)$
preserving the functional $E$. Therefore $E$ projects to a lower
semi-continuous (see Lemma \ref{lsc}) functional $E: H(M)/G \to \R_+$,
where we give $H(M)/G$ the quotient topology. Our main objective is to prove

\begin{thm}\label{T11}
1. $e(M)\ge \mu>0$, where $\mu$ depends only on $R$, $\la_1(M)$ and geometry of $M$.

2. If $M$ admits a cocompact isometric group action, then $E: H(M)/G\to \R_+$ is
proper in the sense that
$$
E^{-1}([0, T])
$$
is compact for every $T\in \R_+$. In particular, $e(M)$ is
attained.
\end{thm}

{\em Sketch of the proof.} With every harmonic function $h=h_\chi\in H(M)$ we associate a
finite set $K_I\subset M$ of the centers $x$ of {\em type 1 special $R$--necks} $N(x)$.
Roughly speaking, these necks encode the partition of $Ends(M)$ into the subsets
$\chi^{-1}(0), \chi^{-1}(1)$. For all but one component $L$ of $M\setminus N(x)$,
$\chi|Ends(L)$ is constant. We also verify that
the function $\chi$ is constant on $Ends(M')$ for each component $M'$ of
$$
M\setminus \bigcup_{x\in K_I} int(N(x)).
$$
Every neck $N(x)$ centered at $x\in K_I$, ``contributes''
at least $\mu>0$ to the energy of $h$. This establishes the inequality
$$
e(M)\ge \mu>0.
$$
If $E(h)\le E$, we also obtain an upper bound on the
cardinality of $K_I$: $|K_I|\le \kappa_1(E)$.

Suppose that $h_n=h_{\chi_n}\in H(M)$ is a sequence of functions with $E(h_n)\le E$.
The corresponding sets $K^{(n)}=K_I(\chi_n)$ break into subsets $K^{(n)}_{i}$ of uniformly
bounded diameter,  so that the distance between distinct subsets diverges to infinity as
$n\to\infty$. Using the group $G$, we normalize the functions $h_n$ so that $K^{(n)}_{1}$ is
contained in a fixed compact subset of $M$. Then the sequence $(h_n)$ subconverges to
a harmonic function $h: M\to [0,1]$. Since each neck $N(x), x\in K^{(n)}_{1}$,
contributes at least $\mu$ to the energy of each function $h_n$, we conclude that
$E(h)\ge \mu>0$. Lastly, we need to check that $h$ extends to a function
$\chi: Ends(M)\to \{0,1\}$ (a priori, this extension
might attain other values in $[0, 1]$ as well).
This follows from the ``uniform connectedness'' considerations
and uniform estimates for the behavior of the functions $h_n$ at the points
far away from $K^{(n)}$.

\medskip
Although it is not needed for the group--theoretic applications, we will
also prove

\begin{thm}\label{T12}
(Finiteness theorem.) Suppose that $M$ admits a cocompact isometric group action.
Then $H(M)$ contains only finitely many $G$--orbits of functions
$h\in H(M)$ for which $E(h)< e(M)+\mu/2$.
\end{thm}

{\bf Acknowledgements.} This paper was motivated by numerous
discussions with Mohan Ramachandran,
to whom I am grateful for many valuable references and suggestions.
During the writing of this paper the I was partially supported by
the NSF grants DMS-04-05180 and DMS-05-54349. This paper
was written when I was visiting
Max Plank Institute for Mathematics in Sciences in Leipzig.

\section{Preliminaries} \label{prelim}

Throughout this paper, we let $M$ be a complete Riemannian manifold of bounded geometry, i.e.,
its injectivity radius is bounded from below by some $C_1>0$ and the
absolute value of the sectional curvature is bounded from above by some
$C_2<\infty$. We say that a constant $C$ depends only on
{\em geometry of $M$} if it depends only on dimension of $M$, and the numbers
$C_1$ and $C_2$.

\begin{notation}
For a subset $N\subset M$ let $N^c$ denote $M\setminus int(N)$.
\end{notation}

\begin{notation}
Given a subset $N\subset M$, let $B_R(N)$ denote the collection of
points in $M$ which are within distance $\le R$ from $N$.
Thus, $B_R(x)$ is the closed $R$-ball centered at $x$.
\end{notation}

\begin{notation}
For subsets $S, T\subset M$, define
$$
\dist(S, T):= \inf \{d(x, y): x\in S, y\in T\}.
$$
\end{notation}

We will assume that $M$ has infinitely many ends. We say that a metric
ball $N=B_r(x)\subset M$ is an $r$-{\em neck} if $N^c$ has at
least three unbounded components.

\begin{assumption}
There exist a number $R$ such that $R$-necks cover $M$.
\end{assumption}

For instance, this assumption holds if $M$ admits a
cocompact isometric group action. We fix $R$ satisfying the above assumption
from now on and will refer to $R$-necks simply as {\em necks}.

\begin{thm}\label{L11}
Under the above assumption, $M$ is non-amenable, i.e., its Cheeger constant is positive:
$$
\eta(M)=\inf \left\{ \frac{Area(\D C)}{Vol(C)}: C\subset M \right\} >0.
$$
Here the infimum is taken over all compact subsets $C\subset M$ with
piecewise-smooth boundary and nonempty interior.
\end{thm}
\proof See \cite{Pittet}. \qed

Let $\la_1(M)$ denote the 1-st eigenvalue of $M$. Then, by Cheeger's theorem
(see \cite[Page 91]{Schoen-Yau}), we have
$$
\la_1(M)\ge \eta^2(M)/4.
$$
In particular,
$$
\la_1(M)>0.
$$

\begin{thm}
Let $M$ be a Riemannian manifold of bounded geometry, so that $\la_1(M)>0$.
Then,  every continuous function $\chi: Ends(M)\to \{0, 1\}$, admits a continuous
extension to a (unique) function
$$
h=h_\chi: \ol{M} \to [0, 1]
$$
whose restriction to $M$ is harmonic.
\end{thm}
\proof This theorem was proven by Kaimanovich and Woess in
\cite[Theorem 5]{Kaimanovich-Woess} using probabilistic methods (they
also proved it for functions $\chi$ with values in $[0, 1]$).
In the the context of K\"ahler manifolds, the theorem was proven in
\cite[Theorem 2.6]{Napier-Ramachandran}. In Section \ref{har}, we will
present a proof of this theorem provided by Mohan Ramachandran. \qed

\medskip
Suppose that $\chi_1, \chi_2: Ends(M)\to \{0, 1\}$ are such that
$$
\chi_1\le \chi_2.
$$
Then, by the maximum principle,
$$
h_{\chi_1}\le h_{\chi_2}.
$$
If the equality is attained at some point of $M$, then $\chi_1=\chi_2$.

\medskip
We now restrict to  continuous functions $\chi: Ends(M) \to \{0, 1\}$.

\begin{lem}\label{fe}
Each function $h=h_\chi$ has finite energy
$$
E(h)= \int_M |\nabla h|^2.
$$
\end{lem}
\proof The assertion follows immediately from
Lemma 5.3 (i) \cite[Page 71]{Schoen-Yau}. \qed 

\medskip
A  subset of $Ends(M)$ is called {\em clopen} if it is both open and closed.

\begin{definition}[Clusters]
A clopen subset of $\chi^{-1}(i)$ is called an $i$--{\em cluster}
with respect to the function $\chi$. When $i$ is irrelevant,
we refer to an $i$--cluster as a {\em cluster}.
\end{definition}

\medskip
A {\em domain} in $M$ is a connected
properly embedded codimension $0$ submanifold
$M'\subset M$, which has smooth compact boundary. Then
$$
Ends(M')\subset Ends(M)
$$
and
$$
\ol{M'}\subset \ol{M}
$$
are continuous embeddings.

\begin{definition}
A domain $M'\subset M$ {\em cobounds} an $i$--cluster
(with respect to the function $\chi$) if $Ends(M')$ is a cluster.
\end{definition}

\section{Uniform connectivity}\label{connectivity}

Fix $R>0$ and let ${\mathcal N}\subset M$ be a $\delta$-separated net ($\delta>0$).
In case when $M$ admits a cocompact isometric action of a discrete group $G$, we assume that
$\n$ is $G$--invariant.
Pick a subset $K\subset {\mathcal N}$ of diameter $\le r$ and consider its
$R$-neighborhood $N=B_R(K)$ in $M$.  Define $\Phi(K, r)$ as follows.
For each component $C$ of $N^c$, consider the induced path--metric on $C$.
Then let $\Phi(K,r)$ be the maximum (over all $C$'s)
of the diameters of $C\cap N$ with respect to this metric. In other words,
$\Phi(K,r)$ equals
$$
\max_C \sup\{ x,y\in C\cap N: \inf_{p\in \Pi_{xy}} \length(p)\}
$$
where the $\Pi_{xy}$ is the set of all paths in $C$ connecting $x$ to $y$.

We define the {\em uniform connectivity function}
$$
\phi(r):= \sup \{\Phi(K,r): K\subset \n, \diam(K)\le r\}.
$$
Then $\phi$ is an increasing function. The following lemma is clear:

\begin{lem}\label{L21}
Suppose that $M$ admits an isometric cocompact group action preserving ${\mathcal N}$.
Then $\phi(r)$ is finite for each $r\in \R$.
\end{lem}

In general, $\phi(r)$ need not be finite.

\begin{example}
Let $R=1$.
Start with the complex plane $\C$ with its flat metric. Let $S_n$ denote the double of
$$
\C\setminus (B_1(n)\cup B_1(-n))
$$
across its boundary. Smooth out this metric along the  boundary of $B_1(n)\cup B_1(-n)$
to make it Riemannian. Then $\phi(4)\ge n$ for $S_n$.

Lastly, take the connected sum of the surfaces $S_n$ ($n\ge 3$) as follows:
Remove from each $S_n$  one copy of the disk $D_n=B_1(0)$ and glue $S_n$ to $S_{n+1}$ along
the boundaries of $D_n, D_{n+1}$. Smooth out the resulting metric.
This infinite connected sum has infinite $\phi(4)$.
\end{example}

\begin{assumption}\label{csc}
From now on we assume that $M$ is such that $\phi(r)$ is finite
for each $r\in \R$.
\end{assumption}

One can easily see that finiteness of $\phi$ is independent of the choice of the net
${\mathcal N}$, number $R$
and is invariant under quasi-isometries. (We do not need these properties.)

\medskip

Let $K^{(n)}=\{x_{n,1},...,x_{k,n}\}$ denote a sequence of subsets of cardinality $\le k$ in $M$.
Since $[0,\infty]^{k^2}$ is compact, after passing  to a subsequence, we can assume
that for each $i, j$, there is a limit
$$
\lim_n d(x_{n,i}, x_{n,j})\in [0,\infty].
$$

Thus, we obtain

\begin{lem}\label{break1}
After passing to a subsequence in the sequence $(K^{(n)})$, we can break $K^{(n)}$
as the disjoint union of nonempty subsets
$$
K^{(n)}=\bigcup_{i=1}^l K^{(n)}_{i}
$$
so that

1. $\diam(K^{(n)}_{i})\le D<\infty$, for all $i=1,...,l$, $n\in \N$.

2. $\lim_{n\to\infty} \dist(K^{(n)}_{i}, K^{(n)}_{j})=\infty$ for $i\ne j$.
\end{lem}

\noindent When $n$ is sufficiently large, we obtain that for all $i\ne j$,
\begin{equation}
\label{br}
\dist(B_R(K^{(n)}_{i}), B_R(K^{(n)}_{j}))> d:=\max_{m} \phi(\diam(K^{(n)}_{m})).
\end{equation}

We now take one of the sets $K=K^{(n)}$ and its partition
$$
K^{(n)}=\bigcup_{i=1}^l K^{(n)}_{i}$$
as in the above lemma. By abusing the notation, we will abbreviate $K^{(n)}_{i}$ as $K_i$,
$i=1,..,l$.

Consider the covering of $M$ by the sets $B_R(K_1)$,..., $B_R(K_l)$ and
by the connected components $C_1$,..., $C_m$ of $B_R(K)^c$. Then the nerve of this covering
is a finite graph $\Ga$ without loop and bigons. We will use the notation $K_i, C_j$ for
the vertices of this graph corresponding to the sets $K_i, C_j$.

We will say that $\Ga$ is {\em dual} to $K$.

\begin{lem}\label{L23}
The graph $\Ga$ is a tree provided that (\ref{br}) holds.
In other words, whenever $x, y\in M$ are disconnected by $int(B_R(K))$,
there exists $i$ so that $x, y$ are disconnected by $int(B_R(K_i))$.
\end{lem}
\proof Suppose that $\Ga$ is not a tree. Then it contains a shortest cycle which we denote
$$
K_1-C_1-K_2-....-C_s-K_1.
$$
Let $x\in C_1\cap B_R(K_1), y\in C_s\cap B_R(K_1)$. Then $x$ and $y$ belong to the same
connected component of $M\setminus int(B_R(K_1))$. Therefore, by the inequality (\ref{br}), there
exists a path $p$ in  $M\setminus int(B_R(K_1))$ disjoint from
$$
B_R(K \setminus K_1)$$
connecting $x$ and $y$. See Figure \ref{f2.fig}.
Therefore this path has to be contained in both $C_1$ and $C_s$.
Hence, $C_1=C_s$. Contradiction. \qed

\begin{figure}[tbh]
\begin{center}
\input{f2.pstex_t}
\end{center}
\caption{\sl   }
\label{f2.fig}
\end{figure}

\section{Estimates on harmonic functions on $M$}

{\bf Gradient estimate}, see \cite[Page 17]{Schoen-Yau}.
There exists a constant $C=C_{grad}$ which depends only on geometry
of $M$, so that for every positive harmonic function $u: M\to \R$ we have
$$
|\nabla u(x)|\le C u(x)
$$
for all $x\in M$.

\bigskip
{\bf Decay estimates for harmonic functions.}

\begin{prop}\label{P31}
There exists a function $\rho(\eps, D, k), \eps>0, D>0$, which depends only on
the  geometry of $M$, so that the following holds.

Let $M'\subset M$ be a domain whose boundary $\D M$
is the union of at most $k$ subsets $\D_i M'$, each of diameter $\le D$.
Set $r:=\rho(\eps, D, k)$. Let $h: M\to (0, 1)$ be a harmonic function which vanishes on
$Ends(M')$. Then:

For every $x\in T:=M'\setminus B_{r}(\D M')$, we have
$$
h(x)\le \eps.
$$

\end{prop}
\proof Given the fact that $\la_1(M)>0$, the proof follows by repeating the
arguments of Lemma 5.3 (part (iii)) in   \cite[Chapter II]{Schoen-Yau}. (This lemma
establishes uniform exponential decay for harmonic functions which
converge to zero at infinity.) See also \cite[Lemmata 1.1, 1.2]{Li-Wang}. \qed

\begin{cor}\label{C31}
Suppose that $M_0, M_1\subset M$ are noncompact disjoint
domains, so that $\diam(\D M_i)\le D$, and
$\chi|Ends(M_i)\equiv i$, $i=0, 1$.
Let $\gamma$ denote a shortest geodesic segment connecting
$\partial M_0$ to $\partial M_1$ and let $\length(\gamma)\le l$.
Then
$$
E(h| B_1(\gamma \cup M_0 \cup M_1))\ge \mu(l, D),
$$
where the function $\mu(l, D)>0$ depends only on the geometry of $M$.
\end{cor}
\proof Take $\eps=1/10$.
By applying Proposition \ref{P31} to the functions $h|M_0$
and $(1-h)|M_1$, we find points $x_i\in M_i$ such that
$d(x_i, \D M_i)=\rho=\rho(\eps, D, 1)$ and
$$
h(x_0)\le \eps, \quad  h(x_1)\ge 1-\eps.
$$
It follows that $d(x_0, x_1)\le 2\rho+l$ and
$$
|h(x_0)-h(x_1)|\ge 1-2\eps.
$$

Let $y_i\in \D M'_i$ denote the end-points of $\ga$. Connect $y_i$ to $x_i$
by the shortest geodesic segments $\al_i, i=0, 1$. Let
$\be:= \al_0\cup \ga \cup \al_1$. Then length of $\be$
is at most $2\rho+l$.

By the mean value theorem, there exists a point $y\in \be$ so that
$$
|\nabla h(y)|\ge \frac{1-2\eps}{2\rho+l}.
$$
Therefore
$$
E(h|B_1(y))\ge \mu(l, D)=Const \frac{0.64}{(2\rho+l)^2},
$$
where $Const$ depends only on geometry of $M$. \qed

\medskip

\begin{lem}\label{lsc}
The energy function $E: H(M)\to \R_+$ is lower semi--conti\-nuous.
\end{lem}
\proof Let $h=h_\chi\in H(M)$ be the limit
$$
h=\lim_{n\to\infty} h_n, \quad h_n\in H(M).
$$
Let $\eps>0$.
Pick a sufficiently large ball $B_r(o)\subset M$, so that each unbounded component $M_i$
($i=1,...,q$) of $B_r(o)^c$ cobounds a cluster with respect to $\chi$.
Then, since $E(h)$ is finite (Lemma \ref{fe}), there exists $\rho\ge r$, so that
for each $i$,
$$
E(h|M_i\setminus B_\rho(o))\le \eps.
$$

Let $C$ denote the compact in $M$ which is the union of $B_\rho(o)$ and the compact
components of $B_r(o)^c$. As uniform convergence $h_n|C$ implies uniform convergence of
these functions  in $C^1$-norm (by the gradient estimate), we obtain
$$
E(h|C)=\lim_{n\to\infty} E(h_n|C).
$$
Therefore,
$$
E(h)\le E(h|C)+ q\eps \le  \lim inf_{n\to\infty} E(h_n).
$$
Since $q$ is constant and $\eps$ is arbitrarily small, we obtain
$$
E(h)\le \lim inf_{n\to\infty} E(h_n).  \qed
$$

\section{An existence theorem for harmonic functions}\label{har}

\begin{thm}
Let $\chi: Ends(M)\to \{0, 1\}$ be a continuous function. Then $\chi$
admits a harmonic extension to $M$.
\end{thm}
\proof (M. Ramachandran.)
Let $\varphi$ denote a smooth extension of $\chi$ to $M$
so that $d\varphi$ is compactly supported.

We let $W^{1,2}_o(M)$ denote the closure of $C_c^\infty(M)$ with respect
to  the norm
$$
\|u\|:= \|u\|_{L_2} + \sqrt{E(u)}.
$$

Consider the affine subspace of functions
$$
\f:= \varphi+ W^{1,2}_o(M)\subset L_{loc}^2(M).
$$
Then the energy is well-defined on $\f$ and we
set $E:= \inf_{f\in \f} E(f)$.

Note that, since $\f$ is affine, for $u, v\in \f$ we also have
$$
\frac{u+v}{2}\in \f,
$$
in particular,
$$
E(\frac{u+v}{2})\ge E
$$
and we set
$$
E(u, v):=  2E(\frac{u+v}{2})- \frac{E(u)+ E(v)}{2}.
$$
The latter equals
$$
E(u, v):= \int_M \< \nabla u, \nabla v\>
$$
in the case when $u, v$ are smooth. We thus obtain
$$
E(u,v)\ge 2E- \frac{E(u)+ E(v)}{2}
$$
for all $u, v\in \f$. Hence,
\begin{equation}
\label{f1}
E(u-v)=E(u) +E(v) -2 E(u, v)\le 2 E(u) + 2 E(v) - 4E.
\end{equation}

Pick a sequence $u_n \in \f$ such that
$$
\lim_{n\to\infty} E(u_n)=E.
$$
Then, according to (\ref{f1}),
$$
E(u_m- u_m) \le 2 E(u_n) + 2 E(u_m) - 4E= 2 (E(u_n)-E) + 2(E(u_m)-E).
$$

Since $\la:=\la_1(M)>0$, we obtain
\begin{equation}
\label{f2}
\la \int_M f^2 \le E(f)
\end{equation}

for all $f\in W^{1,2}_o(M)$. Therefore, the functions
$v_n:= u_n-\varphi\in W^{1,2}_o(M)$ satisfy
$$
\|v_n - v_m\| \le (2+ \la^{-1}) (E(u_n)-E + E(u_m)-E).
$$
Hence, the sequence $(v_n)$ is Cauchy in $W^{1,2}_o(M)$. Set
$$
v:= \lim_n v_n, u:= \varphi+ v\in \f.
$$

By semicontinuity of energy, $E(u)=E$. Therefore, $u$ is harmonic and, hence,
smooth. Since $d\varphi$ is compactly supported, the function $v$ is
also harmonic away from a compact subset $K\subset M$. By the inequality
(\ref{f2}), we have
\begin{equation}
\label{f3}
\int_M v^2 \le \la^{-1} E(v)<\infty.
\end{equation}

Let $r>0$ denote the injectivity radius of $M$. Pick a base-point $o\in M$.
Then (\ref{f3}) implies that
there exists a function $\rho: M\to \R_+$ which converges to $0$ as
$d(x, o)\to \infty$, so that
$$
\int_{B_r(x)} v^2(x) \le \rho(x)
$$
for all $x\in M$. By the gradient estimate, there exists $C_1<\infty$ so that
$$
\sup_{B_r(x)} v^2 \le C_1 \inf_{B_r(x)} v^2
$$
provided that $d(x, K)\ge r$. Therefore,
$$
v^2(x) \le \frac{C_1}{Vol(B_r(x))} \int_{B_r(x)} v^2 \le C_2 \rho(x).
$$
Thus
$$
\lim_{d(x, o)\to\infty} v(x)=0.
$$

Therefore the harmonic function $u$ extends to the function $\chi$ on
$Ends(M)$. \qed

\section{Geometry of necks}\label{necks}

Let $R$ be as in Section \ref{prelim}.
Pick ${\mathcal N}\subset M$, a $\delta$--separated $R$--net in $M$.
If $M$ admits an isometric cocompact action $G\acts M$, we assume that this
net is $G$--invariant. For $x\in \n$ we let $N(x):= B_R(x)$
denote the corresponding neck.

\begin{definition}
Given a nonconstant function $\chi: Ends(M)\to \{0, 1\}$,
we say that a neck $N=N(x)$ is a {\em regular} $\theta$--neck ($\theta\in \{0, 1\}$)
if for all but one components $M'$ of $N^c$ satisfy
$$
\chi|Ends(M')\equiv \theta.
$$

A neck which is not regular, is called {\em special}.
\end{definition}

\begin{figure}[tbh]
\begin{center}
\input{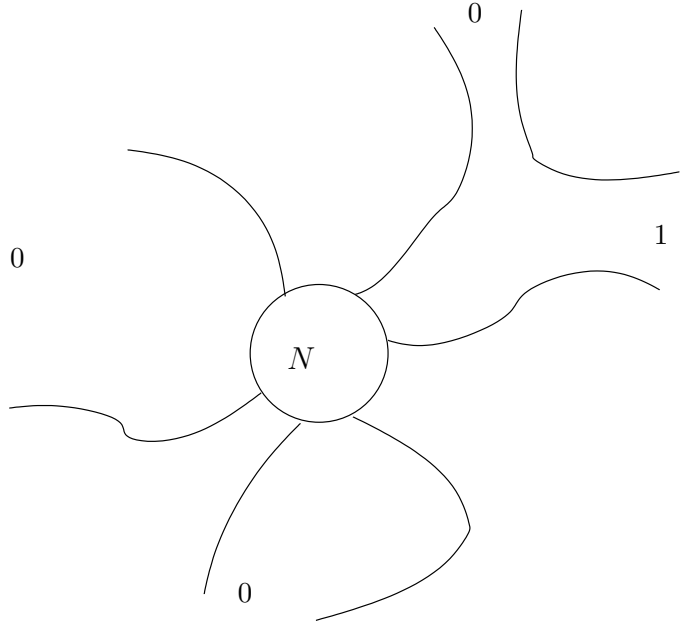}
\end{center}
\caption{\sl Regular neck. }
\label{f3.fig}
\end{figure}

There are two types of special necks:

{\bf Type 1.} There exists at most one (unbounded) component $M'$ of $N^c$
which does not cobound a cluster and there are at least two (unbounded)
components $M_0, M_1$ of $N^c$ so that $\chi|Ends(M_i)\equiv i, i=0, 1$.

{\bf Type 2.} There are at least who components $M_1, M_2$ of $M\setminus int(N)$
which do not cobound clusters.

\begin{figure}[tbh]
\begin{center}
\input{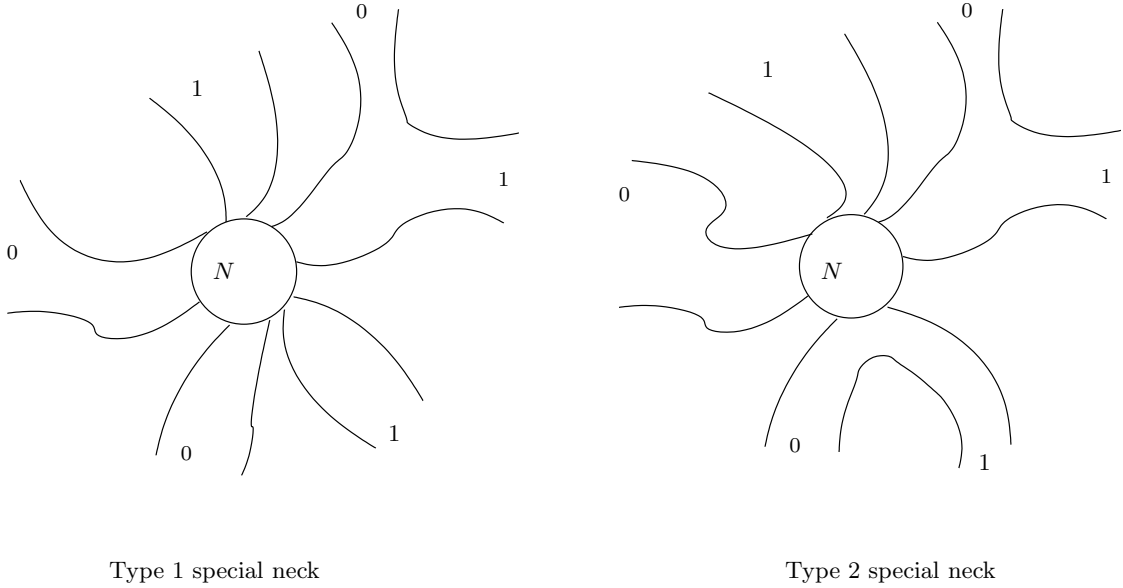}
\end{center}
\caption{\sl Special necks. }
\label{f4.fig}
\end{figure}

\medskip
Let $K\subset \n$ denote the set of centers of special necks and let
$K_{I}$ and  $K_{I\! I}$ denote the subsets of $K$
consisting of the centers of type 1 and 2  necks respectively.

\begin{rem}\label{L51}
Suppose that there exists a special neck $N(x)$ of type 1,
so that each component of $N(x)^c$ cobounds a cluster.
Then every neck $N(y)$ disjoint from $N(x)$ is regular.
\end{rem}

\begin{lem}\label{L52}
Suppose that $N_i=N(x_i)$ are regular $\theta_i$-necks, $i=1, 2$, which have nonempty
intersection. Then $\theta_1=\theta_2$.
\end{lem}
\proof We will consider the most interesting case, when both $N_i^c$ contain
exactly one complementary component $M_i'$ which does not cobound a cluster and will
leave the remaining cases to the reader.

Then
$$
M_2'\cup (N_2\setminus N_1) \subset M_1'.
$$
Suppose that $\theta_1\ne \theta_2$. Let $M_2\subset N_2^c$ be an
(unbounded) component. If $M_2$ is not contained in $M_1'$, then it is contained
in a component $C$ of $N_2^c$ so that $\chi|Ends(C)\equiv\theta_2$. Therefore
$\chi|Ends(C)\equiv \theta_2$ which contradicts our assumption that $\theta_1\ne
\theta_2$.

Hence, we have
$$
M_2\subset M_1'.
$$
Similarly, every component $M_1\subset N_1^c$ is contained in $M_2'$.
But this implies that all unbounded  components of $N_1^c$
are contained in $M_1'$. Therefore $N_1^c$ has only one unbounded
component, i.e., $M_1'$. Contradiction. \qed

\begin{lem}\label{L53}
Suppose that $x_i$ is sufficiently close (in the topology of $\ol{M}$)
to a point $\xi\in \chi^{-1}(\theta)\subset Ends(M)$.
Then the neck $N(x_i)$ is a regular $\theta$--neck.
\end{lem}
\proof Pick a base--point $o\in M$.
Let $U\subset \ol{M}=M\cup Ends(M)$ be an open neighborhood of $\xi$ so that
$\chi|U\equiv \theta$. Then (by the definition of topology on $\ol{M}$)
there exists $r_0$ such that for all $r\ge r_0$, if $C$ is
a component of $M\setminus B_r(o)$ which intersects $U$, then $Ends(C) \subset U$.

Hence, there exists a neighborhood $V\subset U\subset \ol{M}$ of $\xi$,
so that for each $x\in V\cap M$, one of the unbounded
complementary components $C$ of $B_R(x)$
will contain $B_{r_0}(o)$ and the other unbounded components $C'$ will be such that
$Ends(C')\subset U$. Therefore, $\chi|Ends(C')\equiv \theta$. It follows that
$B_R(x)$ is a regular $\theta$--neck. \qed

\begin{cor}\label{C52}
If $\chi$ is nonconstant, then there exists at least one special neck
in $M$.
\end{cor}
\proof Since $\chi$ is nonconstant, the above lemma implies that
$M$ contains at least one regular $i$-neck for $i=0, 1$.
Now the assertion follows Lemma \ref{L52} and connectedness of $M$. \qed

\begin{cor}\label{C53}
The subset $K\subset \n$ of centers of special
necks $N(x)$ is finite.
\end{cor}
\proof The statement follows from compactness of $\ol{M}$ combined with
Lemma \ref{L53}. \qed

\begin{lem}\label{L63}
Let $N=N(x)$ be a type 2 special neck. Then for every component $M'$ of
$N^c$ which does not bound a cluster, $M'\cup N$ contains a type 1
special neck $N(y)$.
\end{lem}
\proof
Let $K_{I\! I}'$ denote the subset of $K_{I\! I}$ consisting of points $y\in \n$ such that
$N(y)\subset M'\cup N$. Let $k(M')$ denote the cardinality of $K_{I\! I}'$. We prove Lemma
using induction on $k(M')$.

1. Suppose $k(M')=1$, i.e., $K_2'$ consists only of $v$. If $M'\cup N$
contains no special necks
besides $N(x)$, then $M'$ cobounds a cluster. This is a contradiction.
Thus, $M'\cup N$ contains  a type 2 special neck.

2. Suppose  the assertion holds whenever $k(M')\le k$. Consider $M'$ with $k(M')=k+1$.
Let $z\in K_2'$; then the neck $N(z)$ is special of type 2.
At least one of the unbounded components $M''$ of $N(z)^c$ (which does not cobound a
cluster in $M$) is contained in $M'$ and is disjoint from $N(x)$. Then $k(M'')\le k$.
Therefore, by the induction assumption, $M''\cup N(z)\subset M'\cup N(x)$ contains a
type 1 special neck $N(y)$. \qed

\begin{lem}\label{L64}
Every unbounded component of $B_R(K_I)^c$ cobounds a cluster.
\end{lem}
\proof If not, then there exists a component $M'\subset B_R(K_I)^c$ which contains
a type 2 special neck $N(x)$, whose non-cluster complementary component $M''$ is entirely
contained in $M'$. Therefore, according to Lemma \ref{L63}, $M''$ contains a type 1
special neck  $N(y)$. However $y\notin K_I$. Contradiction. \qed

\section{Compactness theorem}\label{minimal}

The goal of this section is to prove Theorem \ref{T11}.
Let $\chi: Ends(M)\to \{0, 1\}$ be a nonconstant continuous function.
Let $K=K(\chi)\subset \n$ be as in the previous section.
Define   $\mu=\mu_M:=\mu(2R, R)>0$, where $\mu(\cdot, \cdot)$
is the function defined in Corollary \ref{C31}.

\begin{lem}\label{L61}
For  $x\in K_{I}$, let
$M_i(x)$ denote components of $N(x)^c$ which cobound $i$--clusters, $i=0, 1$.
Then
$$
E(h|M_0(x)\cup M_1(x) \cup N(x))\ge \mu.
$$
\end{lem}
\proof The assertion immediately follows from  Corollary \ref{C31}. \qed

\begin{cor}\label{C61}
If $M'$ is a component of $N(x)^c$ which does not cobound a cluster, then
$$
E(h_\chi|M'\cup N(x))\ge \mu.
$$
\end{cor}
\proof Since $M'$ does not cobound a cluster, there exists a special neck $N(y)$
contained in $M'':=M'\cup N$. If this special neck is of type 1, we are done by
Lemma \ref{L61}. If $N(y)$ is of type 2, then, by Lemma \ref{L63}, $M''$
contains a special neck of type 1. Hence, we are again done by Lemma \ref{L61}.
\qed

\begin{cor}
For every $h\in H(M)$, we have
$$
E(h)\ge \mu.
$$
Thus, $e(M)\ge \mu>0$.
\end{cor}

\begin{lem}\label{L62}
There exists a function $\kappa_1(E)$ (which depends only on geometry of $M$)
such that if $E(h_\chi)\le E$, then the cardinality of $K_{I}$ is at most $\kappa_1(E)$.
\end{lem}
\proof For $x\in K_I$ let $M_i(x)$ denote components
of $N(x)^c$ which cobound $i$--clusters, $i=0, 1$.
It is clear that if $N(x)\cap N(y)=\emptyset$, then
the four sets
$$
M_i(x), M_i(y), i=0, 1
$$
are pairwise disjoint. It follows from  Lemma \ref{L61}, that
$$
E(h|M_0(x)\cup M_1(x) \cup N(x))\ge \mu= \mu(2R, R)
$$
for every $x\in K_I$. Thus, the cardinality of $K_{I}$ is at most $E/\mu$. \qed

One can also bound the number of type 2 necks as well,
provided that $E(h_\chi)$ is sufficiently small:

\begin{prop}\label{P61}
Suppose that $E(h)< E=e(M)+\mu/2$. There exists a function $\kappa_2(E)$ such that
the cardinality of $K_{I\! I}$ is at most $\kappa_2(E)$.
\end{prop}

We do not need this fact and leave it without a proof. The Proposition follows from the
proof of Finiteness Theorem, see \S \ref{finiteness}.
Observe, however, that if $E(h)$ is large comparing to $e(M)$,
then one cannot have a uniform upper bound on the cardinality of $K_{I\! I}$.

We are now ready to prove properness of the function $E: H(M)/G\to \R_+$, assuming that
$G$ is a discrete subgroup of $\Isom(M)$ which acts cocompactly on $M$.

Suppose that $h_n=h_{\chi_n}\in H(M)$ is a sequence of harmonic functions
with uniformly bounded energy $E(h_n)\le E<\infty$.
For each $n$ we define the set
$$
K^{(n)}=K_I(\chi_n)\subset \n
$$
of centers of special necks of type 1.
By Lemma \ref{L62}, the cardinality of each $K^{(n)}$ is at most $k\le \kappa_1(E)$.
We break each $K^{(n)}$ as the union
$$
K^{(n)}=\sqcup_{i=1}^l K^{(n)}_{i}
$$
as in Lemma \ref{break1}, so that
$$
\diam(K^{(n)}_{i})\le D, \quad \forall n, \forall i,
$$
and
$$
\lim_n \dist(K^{(n)}_{i}, K^{(n)}_{j})=\infty
$$
for $i\ne j$. Let $\Ga=\Ga_n$ denote the dual graph for the above partition of $K^{(n)}$.
Since the number of vertices and edges of $\Ga$ is uniformly bounded, after
passing to a subsequence we can assume that $\Ga$ does not depend on $n$.

By applying elements of $G$ and passing to a subsequence,
we can assume that a certain point $x_{1n}\in K^{(n)}_{1}$ is a point $o\in \n$ which
does not depend on $n$. Therefore, without loss of generality, we may assume that
$K^{(n)}_{1}$ does not depend on $n$ either.

Let $M_{n,1},...,M_{n,s}$ denote the unbounded components of $B_R(K^{(n)})^c$ which
are adjacent to $B_R(K^{(n)}_{1})$.
Let $M_1',...,M_t'$ be the unbounded components of $B_R(K^{(n)}_{1})^c$ which
are adjacent to $B_R(K^{(n)}_{1})$.
Because $\Ga$ is a tree, it follows that $s=t$ and that distinct components $M_{n,i}$
lie in distinct components $M_{i}'$ for every sufficiently large $n$, and all
$i=1,...,t$. See Figure \ref{f5.fig}.

\begin{figure}[tbh]
\begin{center}
\input{f5.pstex_t}
\end{center}
\caption{\sl   }
\label{f5.fig}
\end{figure}

Recall that each $M_{n,i}$ cobounds a cluster in $Ends(M)$
(with respect to $\chi_n$). Let $\theta_i$ denote the constant value of $\chi_n$ on
$Ends(M_{n,i})$. (After passing to a subsequence, we may assume that these constants
are independent of $n$.) Note that, since $K^{(n)}_{1}$ is the set of centers of
type 1 special neck, there are $i, j$ so that $\theta_i\ne \theta_j$.

Since the functions $h_n$ take values in $(0,1)$, by the gradient estimate,
the family $(h_n)$ is equicontinuous. Therefore, there exists a limit
$h:= \lim_n h_n$, which is again a harmonic function.

\begin{lem}
For each $i=1,...,s$,
$$
\lim_{d(x, o)\to\infty} h(x)=\theta_i
$$
where $x\in M_i'$.
\end{lem}
\proof Let $\eps>0$. Pick $x\in M_i'\setminus B_{\rho(\eps, D, k)}(\D M_i')$.
Then, for sufficiently large $n\ge n_0$,
$$
x\in B_R(K^{(n)})^c.
$$
Therefore, by Proposition \ref{P31}, for $n\ge n_0$,
$$
|h_n(x)-\theta_i|\le \eps. \qed
$$

Hence, the function $h$ extends to a continuous function $\chi: Ends(M)\to \{0, 1\}$:
$$
\chi|Ends(M_i')\equiv \theta_i.
$$
Since $\theta_i\ne \theta_j$ for some $i, j$, we obtain
that $h\in H(M)$. Since $E$ is lower semicontinuous,
the energy functional $E: H(M)/G\to \R_+$ is proper. It is now clear that $E$
attains the minimum $e(M)>0$. This concludes the proof of Theorem \ref{T11}. \qed

\section{Finiteness theorem}\label{finiteness}

In this section we prove Finiteness Theorem \ref{T12}.
Suppose that there are infinitely many $G$--cosets
of functions $h_n\in H(M)$ with $E(h)<e(M)+\mu/2$. Then, after passing to a
subsequence, and using the notation of the previous section,
$\diam(K^{(n)})\to \infty$ and $h=\lim_n h_n$. As before, we normalize the functions
$h_n$ using the group $G$ and pass to a subsequence, so that
$$
K^{(n)}=\bigcup_{i=1}^l K^{(n)}_{i},
$$
where $K^{(n)}_{1}$ is independent of $n$. Let $h=\lim_n h_n$.
Pick $\eps>0$, so that $\eps<\mu/4$. As in the proof of Theorem \ref{T11},
we get a sufficiently large compact subset $C\subset M$ so that for all $n\ge n_0$ we have
$$
E(h|M\setminus C)\le \eps, \quad E(h_n|M\setminus C)\le \eps,
$$
$$
|E(h|C)-E(h_n|C)|\le \eps.
$$
On the other hand, (for large $n$) $M\setminus C$ contains at least one type 1 special neck
$N:=N(x_{n})$, $x_n\in K^{(n)}\setminus K^{(n)}_{1}$. Let $M_i(x_n)$, $i=0,1$
denote the components of $N^c$ which cobound $i$--clusters with respect to $\chi_n$.
Then
$$
U_n:=N\cup M_0(x_n)\cup M_1(x_n)
$$
is disjoint from the compact $C$. According to Lemma \ref{L61},
$$
E(h_n|U_n)\ge \mu.
$$

Putting these inequalities together, we obtain
$$
E(h)\le E(h|C)+\eps \le E(h_n)+2\eps -\mu<e(M),
$$
since $E(h_n)<e(M)+\mu/2$. However, $h\in H(M)$ and $e(M)=\min \{E(h), h\in H(M)\}$.
Contradiction. \qed

\section{Proof of Stallings' theorem}

The goal of this section is to present the rest of Gromov's proof of the
Stallings' theorem on groups with infinitely many ends. The following was proven
by Stallings \cite{Stallings} for torsion-free groups, his proof was extended
 by Bergman \cite{Bergman} to  groups with torsion:

\begin{thm}[Stallings, Bergman]
Let $G$ be a finitely--generated group with infinitely many ends. Then $G$ splits
nontrivially as a graph of groups with finite edge groups.
\end{thm}
\proof Our argument is a slightly expanded version of Gromov's proof in
\cite[Pages 228--230]{Gromov}. Since $G$ is finitely--generated, it admits a
cocompact isometric properly discontinuous action $G\acts M$
on a connected Riemannian manifold $M$. For instance, if $G$ is $k$--generated,
and $F$ is a Riemann surface of genus $k$, we have an epimorphism
$$
\phi: \pi_1(F)\to G.
$$
Then $G$ acts isometrically and cocompactly on the covering space $M$ of $F$ so that
$\pi_1(M)=\ker(\phi)$. Thus, $M$ has infinitely many ends. The manifold $M$
has bounded geometry since it covers a compact Riemannian manifold.

Let $H(M)$ denote the space of harmonic functions
$h: M \to (0,1)$ as in the Introduction. According to Theorem \ref{T11},
there exists a function $h\in H(M)$ with minimal energy $E(h)=e(M)>0$.
Then, for every $g\in G$, the function
$$
g^*h:= h\circ g
$$
has the same energy as $h$ and equals
$$
h_{g^*(\chi)}.
$$
For $g\in G$, define
$$
g_+(h):= \max(h, g^*(h)), \quad g_-(h):= \min(h, g^*(h)).
$$
Set
$$
\Lambda:=\{ x: h(x)=g^*h(x)\}= \{x: h(x)=h(g(x))\}\subset M.
$$

\begin{lem}\label{L71}
$$
E(g_+(h))+ E(g_-(h))= 2E(h).
$$
\end{lem}
\proof Without loss of generality, we may assume that $h\ne g^*(h)$.
Then the set $\Lambda$ has measure zero (see e.g. \cite{Hardt-Simon} or
\cite{Bar}). Set
$$
M_-:= \{x\in M: h(x)> g^*h(x)\}, M_+:= \{x\in M: h(x)< g^*h(x)\}.
$$
We obtain:
$$
E(g_+(h))+ E(g_-(h))=$$
$$
\int_{M_-} |\nabla h(x)|^2 + \int_{M_+} |\nabla g^*h(x)|^2 +
\int_{M_-} |\nabla g^*h(x)|^2
+\int_{M_+} |\nabla h(x)|^2=$$
$$
 = E(h)+ E(g^*(h))=2E(h). \qed
$$

Note that the functions $g_+(h), g_-(h)$ have continuous extension to
$\ol{M}$ (since $h$ does and $G$ acts on $\ol{M}$ by homeomorphisms).
By construction, the restrictions
$$
\chi_+:= g_+(h)|Ends(M), \quad \chi_-:= g_-(h)|Ends(M)
$$
take the values $0$ and $1$ on $Ends(M)$.
Let
$$
h_\pm:=h_{\chi_\pm}$$
denote the corresponding harmonic functions on $M$. Then
$$
E(h_\pm)\le E(g_\pm(h)),
$$
$$
E(h_+) + E(h_-)\le 2E(h)= 2e(M).
$$

Note that it is, a priori, possible that $\chi_-$ or $\chi_+$ is constant.
Set
$$
G_c:= \{g\in G: \chi_- \hbox{~~or~~} \chi_+ \hbox{~~is constant}\}.
$$

We first analyze the set $G\setminus G_c$. For $g\notin G_c$, both
$h_-$ and $h_+$ belong to $H(M)$ and, hence,
$$
E(h_+)=E(h_-)=E(h)=e(M).
$$
Therefore,
$$
E(g_+(h))=E(h_+), \quad E(g_-(h))=E(h_-).
$$
It follows that $g_\pm(h)$ are both harmonic. Since
$$
g_-(h)\le g_+(h),
$$
the maximum principle implies that either $g_-(h)=g_+(h)$ or
$g_-(h)< g_+(h)$. Hence, the set $\Lambda$ is either empty or equals
the entire $M$, in which case
$g^*(h)=h$. Therefore, for every $g\in G\setminus G_c$ on of the
following holds:

1. $g^*h=h$.

2. $g^*h(x)< h(x), \ \forall x\in M$.

3. $g^*h(x)> h(x), \ \forall x\in M$.

\noindent Thus, the set
$$
L:=h^{-1}\left(\frac{1}{2}\right)
$$
is {\em precisely--invariant} under the elements of $G\setminus G_c$:
for every $g\in G\setminus G_c$, either
$$
g(L)=L
$$
or
$$
g(L)\cap L=\emptyset.
$$

We now consider the elements of $G_c$. Suppose that $g$ is such that $\chi_-=0$. Then
$$
g^*(\chi)\le 1-\chi
$$
and, hence,
$$
g^*(h)\le 1-h.
$$
Since these functions are harmonic, in the case of the equality at some $x\in M$,
by the maximum principle we obtain $g^*(h)=1-h$. The latter  implies
that
$$
g(L)=L.
$$

If
$$
g^*(h)< 1-h
$$
then $g(L)\cap L=\emptyset$. The same argument applies in the case
when $\chi_+$ is constant.

\medskip
To summarize, for every $g\in G$ one of the following holds:
\begin{equation}\label{alt}
g^*h=h, \ \ g^*h< h, \ \ g^*h> h, \ \ g^*h=1-h, \ \ g^*h< 1-h,\ \  g^*h>1-h.
\end{equation}

We conclude that $L$ is precisely--invariant under the action of the
entire group $G$. Moreover,
if $g(L)=L$ then either $g^*h=h$ or $g^*h =1-h$.
Since $L$ is compact, its stabilizer $G_L$ in $G$ is finite.

By construction, the hypersurface $L$ separates $M$
in at least two unbounded components.

Since $L$ is compact, there exists $t\in (0, 1)\setminus \frac{1}{2}$ sufficiently close to
$\frac{1}{2}$, which is a regular value of $h$, so that the hypersurface
$S:= h^{-1}(t)$ is still precisely--invariant under $G$. Let $G_S\subset G_L$ denote
the stabilizer of $S$ in $G$.

It is now rather standard
that $G$ splits nontrivially over a subgroup of $G_S$.
We present a proof for the sake of completeness.
(The proof is straightforward under the assumption that $S$ is connected, but
requires extra work in general.)  We proceed by constructing a simplicial $G$--tree $T$
on which $T$ acts without inversions, with finite edge--stabilizers and without a
global fixed vertex.

\bigskip
{\bf Construction of $T$.} Consider the family of functions $\f=\{f=g^*h: g\in G\}$.
Each function $f\in \f$ defines the {\em wall} $W_f=\{x: f(x)=t\}$
and the {\em half-spaces} $W^+_f:=\{x: f(x)>t\}$, $W^-_f:=\{x: f(x)<t\}$ (these spaces
are not necessarily connected).

Let $\e$ denote the set of walls. We say that a wall $W$ {\em separates} $x, y\in M$ if
$$
x\in W^+_f, \quad y\in W^-_f.
$$

Maximal subsets $V$ of
$$
M^o:=M\setminus \bigcup_{f\in \f} W_f$$
 consisting of points which cannot
be separated from each other by a wall, are called {\em indecomposable} subsets of
$M^o$. Note that such sets need not be connected. Set
$$
\v:=\{ \hbox{indecomposable subsets of} ~~~ M^o\}.
$$
We say that a wall $W$ is {\em adjacent} to $V\in \v$ if $W\cap cl(V)\ne\emptyset$.

The next lemma follows immediately from the inequalities (\ref{alt}), provided that
$t$ is sufficiently close to $\frac{1}{2}$:

\begin{lem}\label{L1}
No wall $W_{f_1}$ separates points of another wall $W_{f_2}$.
\end{lem}

\begin{lem}\label{L2}
1. Let $V\in \v$ and $W\in \e$ be adjacent to $V$. Then,
for each component $C$ of $V$, we have $C\cap W\ne\emptyset$.

2. $W\in \e$ is adjacent to $V\in \v$ iff $W\subset cl(V)$.
\end{lem}
\proof 1. Suppose that $V\subset W^+$.
A generic point $x\in C$ is connected to $W=W_f$ by a gradient curve $p: [0,1]\to M$
of the function $f$. The curve $p$ crosses each wall at most once. Since $V$ is
indecomposable and for sufficiently small $\eps>0$, $p(1-\eps)\in V$, it follows that $p$
does not cross any walls. Therefore the image of $p$ is contained in the closure of $C$
and $p(1)\in W\cap cl(C)$

2. Lemma \ref{L1} implies that for $x, y\in W^+$ (resp. $x, y\in W^-$)
which are sufficiently close to $W$, there is no wall which separates $x$ from $y$.
Therefore, such points $x, y$ belong to the same indecomposable set $V^+$ (resp. $V^-$)
which is adjacent to $W$ and $W\subset cl(V^\pm)$.
Clearly, $V^+, V^-$ are the only indecomposable  sets which are adjacent to $W$. \qed

\medskip
Hence, each wall $W$ is adjacent to exactly two elements of $\v$
(contained in $W^+, W^-$ respectively).
We obtain a graph $T$ with the vertex set $\v$ and edge set $\e$, where a
vertex $V$ is incident to an edge $W$ iff the wall $W$ is adjacent
to the indecomposable set $V$.

From now on, we abbreviate $W_{f_i}$ to $W_i$.

\begin{lem}\label{L3}
$T$ is a tree.
\end{lem}
\proof By construction, every point of $M$ belongs to a wall or to an indecomposable set.
Hence, connectedness of $T$ follows from connectedness of $M$.

Let
$$
W_1 - V_1 - W_2 - .... -W_k - V_k -W_1
$$
be an embedded cycle in $T$. This cycle corresponds to a collection of paths
$p_j: [0,1]\to cl(V_j)$, so that
 $$
 p_j(0)\in W_j, \quad p_{j}(1)\in W_{j+1}, j=1,...,k.
 $$
The points of $p_j([0,1])$ are not separated by any wall, $j=1,...,k$.
By Lemma \ref{L1}, the points $p_j(1), p_{j+1}(0)$ are not separated by
any wall either. Therefore, the
points of
$$
\bigcup_{j=1}^k p_j([0,1])
$$
are not separated by $W_1$. However,
$$
p_1((0,1])\subset W_1^+, \quad p_k([0,1))\subset W_1^-
$$
or vice--versa. Contradiction. \qed

\medskip
We next note that $G$ acts naturally on $T$ since the sets $\f$, $\e$ and $\v$ are
$G$--invariant and $G$ preserves adjacency. If $g(W_f)=W_f$, then $g^*f=f$, which
implies that $g$ preserves $W_f^+, W_f^-$. Hence, $g$ fixes the end-points of the edge
corresponding to $W$, which means that $G$ acts on $T$ without inversions.
The stabilizer of an edge in $T$ corresponding to a wall $W$ is finite, since $W$
is compact and $G$ acts on $M$ properly discontinuously.

Suppose that $G\acts T$ has a fixed vertex. This means that the
corresponding indecomposable subset $V\subset M$ is $G$--invariant. Since $G$ acts
cocompactly on $M$, it follows that $M=B_r(V)$ for some $r\in \R_+$.
The indecomposable subset $V$ is contained in the half-space
$W^+$ for some wall $W$. Since $W$ is compact and
$W^-$ is not, the subset $W^-$ is not contained in $B_r(W)$. Thus $W^- \setminus B_r(V)\ne
\emptyset$. Contradiction.

Therefore $T$ is a nontrivial $G$--tree and we obtain a nontrivial
graph of groups decomposition of $G$ where the edge groups are conjugate
to subgroups of the finite group $G_S$. \qed

\noindent Department of Mathematics, 1 Shields Ave.,\\
University of California, Davis, CA 95616, USA\\
kapovich@math.ucdavis.edu

\end{document}